\documentclass[11pt,reqno]{amsart}
\usepackage{amssymb,graphicx}

\usepackage{amsfonts,amscd,amsthm,amsmath}

\usepackage{mathrsfs}

\usepackage{array}

\newcolumntype{C}[1]{>{\centering\arraybackslash}p{#1}}

\numberwithin{equation}{section}

\begin{document}

\begin{center}

\textbf{\large{New identities for the partial Bell polynomials}}

\vskip 3mm\textbf{Djurdje Cvijovi\'{c}}

\vskip 2mm {\it Atomic Physics Laboratory, Vin\v{c}a Institute of
Nuclear Sciences \\
P.O. Box $522,$ $11001$ Belgrade$,$ Republic of Serbia}\\
\textbf{E-Mail: djurdje@vinca.rs}\\

\vskip 2mm \begin{quotation} \textbf{Abstract.}  A new explicit closed-form formula for the multivariate  $(n, k)$th partial Bell polynomial $B_{n,k} (x_1, x_2, \ldots, x_{n - k + 1})$ is deduced. The formula involves multiple summations and makes it possible, for the first time, to easily evaluate $B_{n,k}$ directly for given  values of $n$ and $k$ ($n\geq k, k =2, 3,\ldots$).  Also, a  new addition formula  (with respect to $k$) is found for the polynomials $B_{n,k}$ and it is shown that they admit a new recurrence relation. Several special cases and consequences are pointed out, and some examples are also given.
\end{quotation}
\end{center}

\vskip 2mm\noindent \textbf{\textit{Key Words and Phrases.}} \begin{Small} Partial Bell polynomial; Recurrence relation; Stirling number of the second kind. \end{Small}

\section{Introduction}

For $n$ and $k$ non-negative integers, the (exponential) $(n, k)$th partial Bell polynomial in the variables  $x_1, x_2,\ldots,x_{n - k + 1}$ denoted by $B_{n,k} \equiv B_{n,k}(x_1, x_2, \ldots, x_{n - k + 1})$ may be defined  by the formal power series expansion  \textup{(}see, for instance, \cite[pp. 133, Eq. (3a')]{Comtet}\textup{)}
\begin{equation} \frac{1}{k!} \left(\sum_{m \,= 1}^{\infty} x_m \,\frac{t^m}{m!}\right)^k = \sum_{n\,= k}^{\infty} B_{n,k}(x_1, x_2, \ldots, x_{n - k + 1})\,\frac{t^n}{n!}\qquad(k\geq0),
\end{equation}
\noindent or, what amounts to the same, by the explicit formula  \cite[p. 96]{Hazewinkel}
\begin{equation}  B_{n,k}  = \sum \frac{n!}{\ell_1!\, \ell_2! \ldots \ell_{n-k+1}!} \left(\frac{x_1}{1!}\right)^{\ell_1} \left(\frac{x_2}{2!}\right)^{\ell_2}\ldots \left(\frac{x_{n-k+1}}{(n-k+1)!}\right)^{\ell_{n-k+1}},
\end{equation}
where (multiple) summation is extended over all partitions of a positive integer number $n$ into exactly $k$ parts (summands), {\em i.e.}, over all solutions in non-negative integers $\ell_{\alpha},$ $1\leq\alpha\leq n-k+1,$ of a system of the two simultaneous equations
\begin{equation*}\ell_1 + 2 \,\ell_2 + \cdots + (n-k+1) \,\ell_{n-k+1} = n
\end{equation*}
\noindent and
\begin{equation*}\ell_1 + \ell_2 + \cdots + \ell_{n-k+1} = k.
\end{equation*}

For fixed $n$ and $k$, $B_{n,k}$ has positive integral coefficients and is a homogenous  and isobaric polynomial in its $(n-k+1)$ variables $x_1, x_2, \ldots, x_{n-k+1}$ of total degree $k$ and total weight $n$, {\em i.e.}, it is a linear combination of monomials $x_1^{\ell_1} x_2^{\ell_2}\ldots x_{n-k+1}^{\ell_{n-k+1}}$ whose partial degrees and weights are constantly given by $\ell_1 + \ell_2 + \ldots + \ell_{n-k+1} = k$ and $\ell_1 + 2 \ell_2 + \ldots + ({n-k+1})\ell_{n-k+1} = n$. For some examples of these polynomials see Section 3.

The partial Bell polynomials are quite general polynomials, they have a number of applications and more details about them can be found in Bell \cite{Bell}, Comtet \cite[pp. 133--137]{Comtet}, Hazewinkel \cite[pp. 95--98]{Hazewinkel}, Charalambides \cite[pp. 412--417]{Charalambides} and Aldrovandi \cite[pp. 151--182]{Aldrovandi}.  However, the following formulae for $B_{n,k}$
\begin{equation} B_{n,k}= \frac{1}{x_1}\cdot\frac{1}{n-k} \sum_{\alpha\,=1}^{n-k} \binom{n}{\alpha}\left[(k+1)-\frac{n+1}{\alpha+1}\right] x_{\alpha+1} B_{n-\alpha,k},
\end{equation}
\begin{equation} B_{n, k_1 + k_2} = \frac{k_1!\, k_2!}{(k_1 + k_2)!} \sum_{\alpha\,=0}^n \binom{n}{\alpha} B_{\alpha, k_1} B_{n - \alpha, k_2}
\end{equation}
\noindent and
\begin{align} B_{n, k + 1}  =  & \frac{1}{(k+1)!} \underbrace{\sum_{\alpha_1\,=  k}^{n-1} \, \sum_{\alpha_2\,=  k-1}^{\alpha_1-1} \cdots  \sum_{\alpha_k\, = 1}^{\alpha_{k-1}-1} }_{k }
\overbrace{\binom{n}{\alpha_1}  \binom{\alpha_1}{\alpha_2} \cdots  \binom{\alpha_{k-1}}{\alpha_k}}^{k}\nonumber
\\
&\cdot x_{n-\alpha_1} x_{\alpha_1 -\alpha_2} \cdots x_{\alpha_{k-1}-\alpha_k} x_{\alpha_k} \qquad(n\geq k+1, k\,=1, 2, \ldots)
\end{align}
\noindent appear not to have been noticed in any work on the subject which we have seen. In this note it is aimed to provide short proofs of these results, show some immediate consequences of them and provide some application examples (see also Section 3).

\section{Proof of the main results}

We begin by showing that the identity (1.3) follows without difficulty from the definition of partial Bell polynomials $B_{n,k}$ by means of the generating relation (1.1), given that the next auxiliary result for powers of series is used.  Consider
\begin{equation}\left(\sum_{n\,= 1}^{\infty} f_{n} x^n\right)^k = \sum_{n\,= k}^{\infty} g_{n} (k)\, x^n.\end{equation}
For a fixed positive integer $k,$ we have that:
\begin{equation}  g_{k}(k) = f_{1}^k, \tag{2.2a}
\end{equation}
\begin{equation}g_n (k) = \frac{1}{(n-k) f_1 } \sum_{\alpha\,=1}^{n-k} \Big[(\alpha+1) (k+1)-(n+1)\Big] f_{\alpha+1} \,g_{n-\alpha}(k)\tag{2.2b}
\end{equation}
$$(n\geq k+1).$$

Indeed, by comparing (2.1) with the definition of $B_{n,k}$ in (1.1) and upon setting $g_n(k)= k! B_{n,k}/n!$ and $f_n= x_{n}/n!,$  we arrive at the proposed formula (1.3) by utilizing (2.2b).

Note that (2.2b) may be found in the literature (see \cite{Gould}) but it is not as widely known (and even less used) as it should be. It is exactly for this reason that we derive it starting from the following more general (and equally little known) recurrence relation involving the series coefficients  $f_n$ and $g_n (k)$ in $\left(\sum_{n\,= 0}^{\infty} f_{n} x^n\right)^k = \sum_{n\,= 0}^{\infty} g_{n} (k)\, x^n.$
\begin{equation} \sum_{\alpha\,= 0}^n \big[\alpha (k + 1) - n\big] f_{\alpha} \, g_{n-\alpha}(k) = 0 \qquad(n\geq0).\tag{2.3}
\end{equation}

First, upon taking logarithms of each side of the equation $g(x) = \big[f(x)\big]^k$  and then differentiating both sides of the result with respect to $x$, we obtain $f (x) g'(x) = k \,f'(x) g (x).$  Next, insert the power series expansions of the various functions in this equation and multiply both sides by $x$, to get
\begin{equation*}\sum_{m\,= 0}^{\infty} f_{m} x^m \cdot \sum_{m\,= 0}^{\infty} m \,g_{m}(k) \,x^m = k \sum_{m\,= 0}^{\infty} m\,f_{m} x^m\cdot \sum_{m\,= 0}^{\infty} g_{m}(k)\, x^m. \tag{2.4}
\end{equation*}

\noindent Now, recall that if $\sum\nolimits_{m\,=0}^{\infty} a_m$ and $\sum\nolimits_{m\,= 0}^{\infty} b_m$ are two series, then their Cauchy product is the series $\sum\nolimits_{n\,= 0}^{\infty} c_n$ where $c_n = \sum\nolimits_{k\,= 0}^{n} a_k b_{n-k}$. This is to say that in  the particular case at hand, by equating the coefficients of a given power of  $x$, say $x^n,$ on both sides of (2.4), we have $\sum\nolimits_{\alpha\,= 0}^{n}  (n -\alpha) f_{\alpha}g_{n-\alpha}(k) = k\, \sum\nolimits_{\alpha\,= 0}^{n} \alpha f_{\alpha} g_{n-\alpha}(k),$ which eventually gives (2.3). The recurrence relation (2.3) is clearly valid for an arbitrary real or complex number $k$ and it can be used  to compute successively as many of the unknown  $g_{m}(k)$ values as desired, in order $g_{0}(k), g_{1}(k),g_{2}(k),\ldots$, if $g_{0}(k)$ is known. The special case of (2.3) solved for $g_n(k),$ for $k$ a positive integer and $f_0 \neq 0,$ appears in  various editions of the standard reference book by Gradshteyn and  Ryzhik (see, for instance, \cite[p. 17, Entry 0.314]{Gradshteyn})

Finally, if we suppose $f_0 = 0$ and $f_1\neq0$ then,  from  $\left(\sum_{n\,= 0}^{\infty} f_{n} x^n\right)^k = \sum_{n\,= 0}^{\infty} g_{n}(k)\, x^n,$ where $k$ is a positive integer, it is obvious that the coefficient $g_{n}(k),$ $n = 0, 1, \ldots, k,$ is only nonzero when  $n = k,$ $g_{k}(k)$ then equals $f_{1}^k$ (see (2.2a)), while $\left(\sum_{n\,= 0}^{\infty} f_{n} x^n\right)^k = \sum_{n\,= 0}^{\infty} g_{n}(k)\, x^n$ reduces to (2.1). Therefore, since $g_{0}(k)= g_{1}(k)= \ldots = g_{k-1}(k)=0$ and $f_0 = 0,$ the recurrence relation (2.3) becomes
\begin{equation*}\sum_{\alpha\,= 1}^{n - k} \big[\alpha (k + 1) - n\big] f_{\alpha} \, g_{n-\alpha}(k) = 0 \qquad(n\geq k),
\end{equation*}
\noindent so that, upon replacing $n$ by $n+1$, putting $\alpha + 1$ for $\alpha$ and solving for $g_n(k)$, we have that the coefficients $g_n(k)$, $n\geq k +1,$ are given by (2.2b) above.

In order to prove (1.4) we shall again resort to the generating relation for $B_{n,k}$ (1.1). Let us  by $[t^n] \phi(t)$ denote the coefficient of $t^n$ in the power series of an arbitrary $\phi(t)$. Put $f(t)= \sum_{m \,= 1}^{\infty} x_m \,\frac{t^m}{m!}$, then by (1.1), we have
\begin{equation*} k_1! B_{n,k_1} = n! \,[t^n] f(t)^{k_1} \qquad(n\geq k_1)\end{equation*}
\noindent and
\begin{equation*} k_2! B_{n,k_2} = n! \,[t^n] f(t)^{k_2}\qquad(n\geq k_2),\end{equation*}
\noindent thus
\begin{align} (k_1 + k_2)! B_{n,k_1 + k_2} & =  n! \,[t^n] f(t)^{k_1 + k_2} = n! \,[t^n] \Big( f(t)^{k_1} \cdot f(t)^{k_2}\Big)\nonumber
\\
& \hskip-20mm =  n! \,\sum_{\alpha\,= 0}^n [t^{\alpha}] f(t)^{k_1}\cdot [t^{n- \alpha}] f(t)^{k_2} =  n!\,\sum_{\alpha\,=0}^n \frac{k_1! B_{\alpha,k_1}}{\alpha!}\frac{k_2! B_{n-\alpha,k_2}}{(n-\alpha)!},\tag{2.5}
\end{align}
$$(n\geq k_1 + k_2)$$
\noindent since $[t^n] \Big(\phi(t) \psi(t)\Big) = \sum_{\alpha\,= 0}^n [t^{\alpha}] \phi(t)\cdot [t^{n-\alpha}] \psi(t)$ (the Cauchy product of two series). We conclude the proof by noting that the required expression (1.4) follows by rewriting (2.5).

Lastly, we shall prove the closed-form formula (1.5) by making use of (1.4).  It suffices to show that the addition formula for $B_{n,k}$ (1.4) may be used to deduce the following:
\begin{equation} B_{n,2} = \frac{1}{2!} \sum_{\alpha\,= 1}^{n-1} \binom{n}{\alpha} x_{n - \alpha}\, x_{\alpha} \qquad(n\geq 2),\tag{2.6}
\end{equation}
\begin{equation} B_{n,3} = \frac{1}{3!} \sum_{\alpha\,= 2}^{n-1} \sum_{\beta\,= 1}^{\alpha -1} \binom{n}{\alpha}  \binom{\alpha}{\beta}  x_{n - \alpha} \, x_{\alpha - \beta} \, x_{\beta}\qquad(n\geq 3)\tag{2.7}
\end{equation}
\noindent and
\begin{equation} B_{n,4} = \frac{1}{4!} \sum_{\alpha\,= 3}^{n-1} \sum_{\beta\,= 2}^{\alpha -1} \sum_{\gamma\,= 1}^{\beta-1}\binom{n}{\alpha}  \binom{\alpha}{\beta} \binom{\beta}{\gamma} x_{n - \alpha} \, x_{\alpha - \beta} \,  x_{\beta - \gamma}\, x_{\gamma}\qquad(n\geq 4).\tag{2.8}
\end{equation}

By  bearing in mind that $B_{n,1} = x_{n}$ (this is a simple consequence of the definition $B_{n,k}$ in (1.1))  and upon noticing that $x_{0}= 0$ (again, see (1.1)), the expression for $B_{n,2}$ given in (2.6) follows by (1.4) with $k_1 = 1$ and $k_2 =1$. Further, this result for $B_{n,2}$ together with (1.4), where $k_1 = 2$ and $k_2 =1,$ leads to (2.7). It is clear that by repeating this procedure recursively we may obtain  $B_{n,4},$ and so on.

\section{Further results and concluding remarks}

We remark that  the explicit closed-form formula  for $B_{n,k}(x_1, x_2, \ldots, x_{n - k + 1})$ given by (1.5) is  particularly useful. Namely, it is  hard  to work with the formula (1.2) which explicitly defines $B_{n,k}$ due to complicated multiple summations, and, for instance, it is virtually impossible by its use to write down a polynomial for given values of $n$ and $k$. However, formula (1.5), although  also involves multiple summations, makes this possible. In other words, it is now possible to directly evaluate $B_{n,k}$ for given $n$ and $k$ ($n\geq k, k =2, 3,\ldots$) by utilizing (1.5) instead of computing it  recursively by  making use of some recurrence relations (see, for instance, (1.3)). It is noteworthy to mention that the practical evaluation is greatly facilitated by wide availability of various symbolic algebra programs. In order to demonstrate an application of this result, we list
several of the polynomials $B_{n,k}$ determined by the formula (1.5), where all the computations were carried out by using  Mathematica 6.0 (Wolfram Research)
\begin{small}
\begin{align*}&B_{8,7} = 28 x_1^6 x_2,\qquad B_{9,7} = 378 x_1^5 x_2^2 + 84 x_1^6 x_3,
\\
& B_{10,7} = 3150 x_1^4 x_2^3 + 2520 x_1^5 x_2 x_3 + 210 x_1^6 x_4,
\\
& B_{11,7} = 17325 x_1^3 x_2^4 + 34650 x_1^4 x_2^2 x_3 + 4620 x_1^5 x_3^2 +
 6930 x_1^5 x_2 x_4 + 462 x_1^6 x_5,
\\
& B_{12,7} = 62370 x_1^2 x_2^5 + 277200 x_1^3 x_2^3 x_3 +
 138600 x_1^4 x_2 x_3^2 + 103950 x_1^4 x_2^2 x_4
 \\
 & \quad \quad \,\,+
 27720 x_1^5 x_3 x_4 + 16632 x_1^5 x_2 x_5 + 924 x_1^6 x_6,
 \\
 & B_{13,7} = 135135 x_1 x_2^6 + 1351350 x_1^2 x_2^4 x_3 +
 1801800 x_1^3 x_2^2 x_3^2 + 200200 x_1^4 x_3^3
 \\
 & \quad \quad \,\,+ 900900 x_1^3 x_2^3 x_4 + 900900 x_1^4 x_2 x_3 x_4 +
 45045 x_1^5 x_4^2 + 270270 x_1^4 x_2^2 x_5
 \\
 & \quad \quad \,\, +72072 x_1^5 x_3 x_5 + 36036 x_1^5 x_2 x_6 + 1716 x_1^6 x_7.
\end{align*}
\end{small}

\noindent It should be noted that our results for $B_{8,7}$ $B_{9,7}$  and $B_{10,7}$ are in full agreement with those recorded in the work (for instance) of Charalambides \cite[p. 417]{Charalambides}.

One further  illustration of an application of (1.5) is  the following (presumably) new explicit formula
\begin{align} S(n,k+1)  =  & \frac{1}{(k+1)!} \underbrace{\sum_{\alpha_1\,=  k}^{n-1} \, \sum_{\alpha_2\,=  k-1}^{\alpha_1-1} \cdots  \sum_{\alpha_k\, = 1}^{\alpha_{k-1}-1} }_{k }
\overbrace{\binom{n}{\alpha_1}  \binom{\alpha_1}{\alpha_2} \cdots  \binom{\alpha_{k-1}}{\alpha_k}}^{k}
\end{align}
$$(n\geq k+1, k\,=1, 2, \ldots)$$
\noindent  for the Stirling numbers of the second kind $S(n,k)$ defined by means of (see \cite[Chapter 5]{Comtet}) \begin{equation}S(n,k)=\frac{1}{k!} \sum_{\alpha\,=0}^k (-1)^{k-\alpha} \binom{k}{\alpha}\alpha^n,\end{equation}
\noindent which is an immediate consequence of the  relationship $S(n,k)= B_{n,k}(1,\ldots,1)$ \cite[p. 135, Eq. (3g)]{Comtet}. Moreover, for given $k$, it is easy to sum the multiple sum (3.1) by repeated  use  of the familiar result $(1+x)^n = \sum_{k\,=0}^n \binom{n}{k} x^k$, so that we have:
\begin{align*}
&S(n,2)= \frac{1}{2} \Big(2^n -2\Big) = 2^{n-1}-1,
\\
&S(n,3)= \frac{1}{6} \Big(3^n - 3\cdot 2^n+3\Big),
\\
&S(n,4)=  \frac{1}{24} \Big(4^n - 4\cdot 3^n  +3\cdot 2^{n+1}-4\Big),
\\
&S(n,5)= \frac{1}{120} \Big(5^n - 5 \cdot4^n + 10\cdot 3^n - 10\cdot 2^{n } +5 \Big),
\\
&S(n,6)=  \frac{1}{720} \Big(6^n -6\cdot 5^n + 15\cdot 4^n -20\cdot 3^n + 15\cdot 2^n   -6  \Big),
\end{align*}
\noindent and these expressions agree fully with those which are obtained by using the defining relation (3.2).

\section*{Acknowledgements}\noindent {\small   The author acknowledges financial support from Ministry of Science of the Republic of Serbia under Research Projects 144004 and 142025.}

\end{document}